\documentclass[11pt]{article}
\usepackage{amssymb,amsfonts,amsmath,latexsym,epsf,tikz,url}

\newtheorem{theorem}{Theorem}[section]

\newtheorem{conjecture}[theorem]{Conjecture}
\newtheorem{corollary}[theorem]{Corollary}

\newtheorem{remark}[theorem]{Remark}

\newcommand{\proof}{\noindent{\bf Proof. }}
\newcommand{\qed}{\hfill $\square$\medskip}

\begin{document}

\title{Dominated chromatic number of some operations on a graph}

\author{
Saeid Alikhani$^{}$\footnote{Corresponding author}
\and
Mohammad R. Piri
}

\date{\today}

\maketitle

\begin{center}
Department of Mathematics, Yazd University, 89195-741, Yazd, Iran\\
{\tt alikhani@yazd.ac.ir, piri429@gmail.com}
\end{center}


\begin{abstract}
Let $G$ be a simple graph. The dominated coloring of a graph $G$ is a proper coloring of $G$ such that each color class is dominated by at least one vertex. The minimum number of colors needed for a  dominated coloring of $G$ is called the dominated chromatic number of $G$, denoted by $\chi_{dom}(G)$. In this paper,
we examine the effects on $\chi_{dom}(G)$
 when $G$ is modified by operations on vertex and
edge of $G$.
\end{abstract}


\medskip
\noindent{\bf AMS Subj.\ Class.}: 05C25 

\section{Introduction and definitions}
Let $G=(V, E)$ be a simple graph and $\lambda \in \mathbb{N}$. A mapping $f: V \longrightarrow \{1, 2,... , \lambda\}$ is called a $\lambda$-proper coloring of $G$ if $f(u) \neq f(v)$, whenever the vertices $u$ and $v$ are adjacent in $G$. A color class of this  coloring is a set consisting of all those vertices assigned the same color. If $f$ is a proper coloring of $G$ with the coloring classes $V_1, V_2, ... ,V_{\lambda}$ such that every vertex in $V_i$ has color $i$, sometimes write simply $f= (V_1, V_2, ... ,V_{\lambda})$. The chromatic number $\chi (G)$ of $G$ is the minimum of colors needed in a proper coloring of a graph. 

A dominator coloring of $G$ is a proper coloring of $G$ such that  every vertex of $G$ is adjacent to all vertices of at least one color class. The dominator chromatic number $\chi_{d}(G)$ of $G$ is the minimum number of color classes in a dominator coloring of $G$. The concept of dominator coloring was introduced and studied by Gera, Horton and Rasmussen \cite{e}.
 Let $G$ be a graph with no isolated vertex, the total dominator coloring is a proper coloring of $G$ in which each vertex of the graph is adjacent to every vertex of some (other) color class. The total dominator chromatic number, abbreviated TD-chromatic number, $\chi_{d}^{t}(G)$ of $G$ is the minimum number of color classes in a TD-coloring of $G$. For more information see \cite{nima1,nima2}. 

Dominated coloring of a graph is a proper coloring in which each color class is  dominated by a vertex. The least number of colors needed for a dominated coloring of $G$ is called the dominated chromatic number of $G$ and denoted by $\chi_{dom}(G)$ (\cite{Choopani,cc}). We call this coloring a dom-coloring, simplicity. In the study
of dom-chromatic number of graphs, this naturally raises the question: What happens
to the dom-chromatic number, when we consider some operations on the vertices and
the edges of a graph? In this paper we would like to answer to this question.

\medskip 
In the next section, examine the effects on $\chi_{dom}(G)$
 when $G$ is modified by deleting
a vertex or deleting an edge. In Section $3$, we study the effects on $\chi_{dom}(G)$, when $G$ is
modified by contracting a vertex and contracting an edge. Also we consider another
obtained graph by operation on a vertex $v$ denoted by $G\odot v$ which is a graph obtained
from $G$ by the removal of all edges between any pair of neighbors of $v$ in Section $3$ and
study $\chi_{dom}(G\odot v)$. In the last section we study the dominated chromatic number of subdivision of graph $G$.

\medskip

\section{Vertex and edge removal}

\begin{figure}
	\begin{center}
		\includegraphics[width=0.8\textwidth]{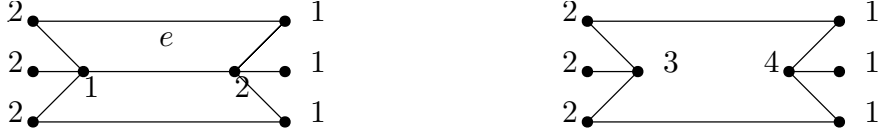}
		\caption{\label{fv} Dominated coloring of $G$ and $G-e$.}
	\end{center}
\end{figure}

The graph $G-v$ is a graph that is made by deleting the vertex $v$ and all edges connected
to $v$ from the graph $G$ and the graph $G-e$ is a graph that obtained from $G$ by simply
removing the edge $e$. We obtain  a bound for
dom-chromatic number of $G-e$ and $G-v$.
\begin{theorem}\label{removeedge}
Let $G$  be a connected graph and $e=uv \in E(G)$ is not a bridge of $G$. Then we have:
\begin{equation*}
\chi_{dom}(G)-1\leq \chi_{dom}(G-e) \leq \chi_{dom}(G)+2.
\end{equation*}
\end{theorem}
\proof 
First we prove $\chi_{dom}(G-e) \leq \chi_{dom}(G)+2$. Suppose that the vertex $v$ has color $i$ and the vertex $u$ has color $j$. We have the following cases: 
\begin{enumerate}
\item[Case 1)]
If $v$ is  the only vertex that dominate the color $j$, with removing the edge $e$ the color class $j$ does not dominate by any vertex, so  in this case we add a new color to the vertex $u$. Therefore $\chi_{dom}(G-e) \leq \chi_{dom}(G)+1$.
\item[Case 2)]
If $u$ is  the only vertex that dominate the color class $j$, then the proof is similar to Case 1.
\item[Case 3)]
If $v$ is the only vertex that dominate color class $j$ and $u$ is the  only vertex that dominate color class $i$, then with removing  edge $e$, we need two new color for the vertices $u$ and $v$. Therefore in this case $\chi_{dom}(G-e) \leq \chi_{dom}(G)+2$.
 \end{enumerate}
Now we prove the left inequality. We shall present a dominated coloring for $G-e$. If we add the edge $e$ to $G-e$, then we have two cases. If two vertices $v$ and $u$ have the same color in the dominated coloring of $G-e$, then in this case, we add a new color, say color $t$  to one of them. Since every color class of $G-e$ dominated by old dominated coloring and color class $t$ dominated by the adjacent vertex, so we have $\chi_{dom}(G)\leq \chi_{dom}(G-e)+1$.
 \qed
 \\
 \begin{remark}
 The lower bound of $\chi_{dom}(G-e)$ in Theorem \ref{removeedge} is sharp. It suffices to consider complete graph $K_3$ as $G$. Also the upper bound is sharp, because as we see in Figure \ref{fv}, $\chi_{dom}(G)=2$ and $\chi_{dom}(G-e)=4$.
 \end{remark} 
 
\medskip
Now we consider the graph $G-v$, and present a lower bound and an upper bound for the dom-chromatic number of $G-v$.
\begin{theorem}\label{removevertex}
Let $G$ be a connected graph, and $v\in V(G)$ is not a cut vertex of $G$. Then we have:
\begin{equation*}
\chi_{dom}(G)-1\leq \chi_{dom}(G-v)\leq \chi_{dom}(G)+deg v -1.
\end{equation*}
\end{theorem}
\proof
First we prove $\chi_{dom}(G)-1\leq \chi_{dom}(G-v)$. We shall present a dom-coloring for $G-v$. If we add a vertex $v$ and all the corresponding edges to $G-v$, then it suffices to give a  new color say $i$ to the vertex $v$. Because  every color class except $i$ dominated by vertices of $G-v$ and color $i$ is  dominated by an adjacent vertex of $v$. This is a dominated coloring for $G$ and therefore we have $\chi_{dom}(G)\leq \chi_{dom}(G-v)+1$.

Now we prove $\chi_{dom}(G-v)\leq \chi_{dom}(G)+deg v -1$. First we give a dom-coloring to $G$. Suppose that the vertex $v$ has the color $i$. We have two cases. If $v$ is the only vertex that dominate all color classes and there is no other vertex with color $i$, then by  removing the vertex $v$, the color classes adjacent with $v$ do not dominate. So we give the new colors $i, a_1, a_2, ... , a_{deg v -1}$  to all adjacent vertices of $v$. Obviously, this is a dom-coloring for $G-v$. So $\chi_{dom}(G-v)\leq \chi_{dom}(G)+deg v -1$. Otherwise, If $v$ is  not only vertex that dominate a color class, then dom-coloring of $G$ is a dom-coloring for $G-v$. So $\chi_{dom}(G-v)\leq \chi_{dom}(G)$. Therefore $\chi_{dom}(G-v)\leq \chi_{dom}(G)+deg v -1$.\qed

\begin{remark} 
	The lower bound in Theorem \ref{removevertex} is sharp. Consider the complete graph $K_n$, as $G$.
	\end{remark} 
	We need the following easy result: 
\begin{theorem}{\rm\cite{piri}}\label{wheel}
If $\Delta (G)=|V(G)|-1$, then $\chi_{dom}(G)=\chi_d^t(G)=\chi (G)$.
\end{theorem}

\begin{theorem}
There is a connected graph $G$, and a vertex $v\in V(G)$ which is not cut vertex of $G$ such that $\vert \chi_{dom}(G)-\chi_{dom}(G-v) \vert$ can be arbitrary large.
\end{theorem}
\proof 
Consider the wheel graph $W_n$. By Theorem \ref{wheel}, 	$\chi_{dom}(W_n)=\chi(W_n)$. We remove the center vertex $v$ of $W_n$. Then $\vert \chi_{dom}(G)-\chi_{dom}(G-v) \vert = \vert \chi_{dom}(W_n)-\chi_{dom}(C_n)\vert$ which can be arbitrary large.\qed

\section{Vertex and edge contraction}
Let $v$ be a vertex in graph $G$. The contraction of $v$ in $G$ denoted by $G/v$ is the graph obtained by deleting $v$ and putting a clique on the (open) neighbourhood of $v$. Note
that this operation does not create parallel edges; if two neighbours of $v$ are already
adjacent, then they remain simply adjacent (see \cite{ars,walsh}). In a graph $G$, contraction of an
edge $e$ with endpoints $u, v$ is the replacement of $u$ and $v$ with a single vertex such that
edges incident to the new vertex are the edges other than $e$ that were incident with $u$
or $v$. The resulting graph $G/e$ has one less edge than $G$. We denote this graph by
$G/e$. In this section we examine the effects on $\chi_{dom}(G)$ when $G$ is modified by an edge
contraction and vertex contraction. First we consider edge contraction:
\begin{theorem}\label{edgecontractin}
Let $G$ be a connected graph and $e\in E(G)$. Then we have:
\begin{equation*}
\chi_{dom}(G)-2\leq \chi_{dom}(G/e) \leq \chi_{dom}(G)+1.
\end{equation*}
\end{theorem}
\proof
First we consider a dom-coloring, say $c$ for $G$. Suppose that the end point of $e$ are the vertices $u$, $v$, the vertex $u$ has the color $i$ and the vertex $v$ has color $j$. We give all used colors in the coloring $c$ to the vertices $V(G)- \lbrace u, v\rbrace$. Now we give the new color $t$ to $u$ and $v$. So every color class that dominated by vertices $u$ and $v$, is dominated by $u$ and $v$, and the color class $t$ dominated by an adjacent vertex of $u=v$ and other color classes dominated by old dom-coloring of $G$. Then this is a dom-coloring $G/e$. If any other vertex does not have the color $i$ and $j$, then it suffices to give color $i$ to one of the adjacent vertices of $u$ (or $v$) in $G$. Then this is a dom-coloring for $G/e$. So we have $\chi_{dom}(G/e) \leq \chi_{dom}(G)+1$.

To prove  the lower bound for $\chi_{dom}(G/e)$, we  give a dom-coloring to $G/e$. We add the removed vertex and all the corresponding edges to $G/e$ and keep the old coloring for the new graph. Now we consider the endpoints of $e$ and remove the used color. Now add new color $i$ and $j$ to these vertices.
Let the vertex $u$ has the color $i$ and the vertex $v$ has the color $j$. So the color class $i$ dominated by $v$ and the color class $j$ dominated by $u$, and all color classes in $V(G)-\lbrace u, v\rbrace$ dominated by dom-coloring of $G/e$. So this is a dom-coloring and we have $\chi_{dom}(G)\leq \chi_{dom}(G/e)+2$. Therefore $\chi_{dom}(G)-2\leq \chi_{dom}(G/e)$.\qed

\begin{remark} 
  By  considering the cycle $C_4$ as $G$ in Theorem \ref{edgecontractin}, we observe that the upper bound  is sharp.
 \end{remark} 
 
 Checking graphs with small order, we have observed that the upper bound 
 in Theorem \ref{edgecontractin} can be decreased by one, but we are not able to prove it. So we state the following conjecture:  
 \begin{conjecture} 
  $\chi_{dom}(G)-1\leq \chi_{dom}(G/e)$.
  \end{conjecture} 
 \begin{corollary}
 Suppose that $G$ is a connected graph and $e\in E(G)$ is not a bridge of $G$. We have
 \begin{equation*}
 \dfrac{\chi_{dom}(G-e)+\chi_{dom}(G/e)-3}{2}\leq \chi_{dom}(G)\leq \dfrac{\chi_{dom}(G-e)+\chi_{dom}(G/e)+3}{2}.
 \end{equation*}
\end{corollary}
\proof
It follows from Theorems \ref{removeedge} and \ref{edgecontractin}.\qed 

\medskip
Now we consider the vertex contraction of graph $G$ and examine the effect on $\chi_{dom}(G)$ when $G$ is modified by this operation.
\medskip
\begin{theorem}\label{vertexcontractin}
Let $G$ be a connected graph and $v\in V(G)$. Then we have:
\begin{equation*}
\chi_{dom}(G)-1\leq \chi_{dom}(G/v)\leq \chi_{dom}(G)+deg v-1.
\end{equation*}
\end{theorem}
\proof
First we give a dom-coloring, say $c$ to $G/v$. We add the vertex $v$, add all the removed edges and remove all the added edges. It suffices to give the vertex $v$ a  new color $i$. All the vertices except the vertex $v$ can use the previous colors in coloring $c$. All color classes which are adjacent to $v$ dominated by $v$ and other color classes dominated by dom-coloring $G/v$. So we have $\chi_{dom}(G)\leq \chi_{dom}(G/v)+1$.\\
To prove  the upper bound, we present a dom-coloring for $G$. We remove the vertex $v$ and create $G/v$. We consider one of the adjacent vertices of $v$ and do not change its color and give the new colors $i, i+1, ... , i+deg (v) -1$ to other adjacent vertices of $v$. All the color classes which are not adjacent to $v$ dominated by previous vertices of $G/v$. Other color classes is dominated by previous vertices that were adjacent with $v$ in graph $G$. So we have $\chi_{dom}(G/v)\leq \chi_{dom}(G)+deg v-1$. Therefore we have the result.\qed
\\
\medskip
 {\bf Remark 4.}
 The bounds in Theorem \ref{vertexcontractin} are sharp. For the upper bound consider the graph $C_4$ and for the lower bound consider $C_5$.
 \medskip 
 \begin{corollary}
 Let $G $ be a connected graph. For every $v\in V(G)$ which is not cut vertex of $G$, we have:
 \begin{equation*}
 \dfrac{\chi_{dom}(G/v)+\chi_{dom}(G-v)}{2}-deg v +1\leq \chi_{dom}(G) \leq \dfrac{\chi_{dom}(G/v)+\chi_{dom}(G-v)}{2}+1.
 \end{equation*}
 \end{corollary}
\proof
It follows from Theorem \ref{removevertex} and \ref{vertexcontractin}.\qed

\medskip
Here we consider another operation on vertex of a graph $G$ and examine the effects
on $\chi_{dom}(G)$ when we do this operation. We denote by $G \odot v$ the graph obtained from $G$
by the removal of all edges between any pair of neighbors of $v$, note $v$ is not removed
from the graph \cite{ars}. The following theorem gives upper bound and lower bound for $\chi_{dom}(G\odot v)$.
\begin{theorem}\label{ope}
Let $G$ be a connected graph and $v\in V(G)$. Then we have:
\begin{equation*}
\chi_{dom}(G)-deg v +1\leq \chi_{dom}(G \odot v) \leq \chi_{dom}(G)+1.
\end{equation*}
\end{theorem}
\proof
First we prove $\chi_{dom}(G)-deg v +1\leq \chi_{dom}(G \odot v)$. Consider the graph $G \odot v$ and shall find a dom-coloring for it. We make $G$ from $G \odot v$ and just change the color of all the adjacent vertices of $v$ except one of them, say  $u$ to the new colors $a_1, a_2, ... , a_{deg v -1}$ and do not change the color of $v, u$ and  other vertices. This is a dom-coloring for $G$, because the vertex $v$ is dominated by an adjacent vertex of $v$, say  $u$ and the color classes adjacent to $v$ is dominated by $v$ and the other color classes is dominated by previous vertices. So we have $\chi_{dom}(G) \leq \chi_{dom}(G\odot v)+deg v -1$.\\
Now we prove $\chi_{dom}(G \odot v) \leq \chi_{dom}(G)+1$. We give a dom-coloring for the graph $G$. Suppose that the vertex $v$ has the color $i$. We have the following cases:
\begin{enumerate}
\item[Case 1)] 
The color $i$ uses only for the vertex $v$. In this case, if a color class dominated by only one of the adjacent vertices of the vertex $v$, then with operation $G\odot v$, is not dominated. In this case we give the new color $j$ to adjacent vertices of vertex $v$. So this is a dom-coloring for $G\odot v$. Therefore $\chi_{dom}(G \odot v) \leq \chi_{dom}(G)+1$.
\item[Case 2)]
The color $i$ uses for another vertex except the vertex $v$. In this case, we give the new color $j$ to all of these vertices (except $v$). This is a dom-coloring for $G\odot v$. So we have $\chi_{dom}(G \odot v) \leq \chi_{dom}(G)+1$ and we have the result.
\end{enumerate}
\qed

\begin{remark} \label{4}
The bounds in Theorem \ref{ope} are sharp. Consider the graph $K_n$ as $G (n\geq 3)$, $\chi_{dom}(K_n)=n$. Now for every $v\in V(K_n)$, $K_n \odot v$ is the star graph $S_n$ and we have $\chi_{dom}(S_n)=2$.
\end{remark} 
 By Remark \ref{4}  we have the following result:
\begin{corollary}
There is a connected graph $G$ and $v\in V(G)$ such that $\dfrac{\chi_{dom}(G)}{\chi_{dom}(G\odot v)}$ can be arbitrary large.
\end{corollary}

\section{Subdivision}
The $k$-subdivision of $G$, denoted by $G^{\frac{1}{k}}$, is constructed by replacing each edge $v_iv_j$
of $G$ with a path of length $k$, say $P^{\{v_i ,v_j\}}$. These $k$-paths are called superedges, any new
vertex is an internal vertex, and is denoted by $x_l^{\{v_i,v_j\}}$
 if it belongs to the superedge
$P^{\{v_i ,v_j\}}, i < j$ with distance $l$ from the vertex $v_i$, where $l\in \{1, 2,..., k-1\}$. Note that
for $k = 1$, we have $G^{\frac{1}{1}} = G^1 = G$, and if the graph $G$ has $n$ vertices and $m$ edges, then
the graph $G^{\frac{1}{k}}$ has $n + (k - 1)m$ vertices and $km$ edges. We need the following theorem: 
\begin{theorem}\rm(\cite{cc})\label{path}
For $n\geq 3$,
 \[
 	\chi_{dom}(P_n)=\chi_{dom}(C_n)=\left\{
  	\begin{array}{ll}
  	{\displaystyle
  		\dfrac{n}{2}}&
  		\quad\mbox{if $n\equiv 0 ~(mod\ 4)$, }\\[15pt]
  		{\displaystyle
  			\lfloor \dfrac{n}{2}\rfloor +1} &
  			\quad\mbox{otherwise. }
  				  				\end{array}
  					\right.	
  					\]
\end{theorem}
 \begin{theorem}\label{frac}
 If $G$ is a connected graph with $m$ edges and $k \geq 2$, then
\begin{equation*}
\chi_{dom}(P_{k+1})\leq \chi_{dom}(G^{\frac{1}{k}})\leq (m-1)\chi_{dom}(P_k)+\chi_{dom}(P_{k+1}).
\end{equation*}
 \end{theorem}
 \proof
 For the right inequality, let $e = uu_1$ be an arbitrary edge of $G$. This edge is
replaced with the superedge $P^{\{u,u_1\}}$ in $G^{\frac{1}{k}}$, with vertices $\{u, x_1^{\{u,u_1\}},..., x_{k-1}^{\{u,u_1\}}, u_1 \}$.
We color this superedge with $\chi_{dom}(P_{k+1})$ colors as dominated coloring of $P_{k+1}$
(Theorem \ref{path}). If $N_G(u) = \{u_1,..., u_s\}$ then we color the vertices of paths $P^{\{u,u_i\}}$
such that
the color of $u$ in $P^{\{u,u_i\}}$, for any $2 \leq i \leq s$, is the same as color $u$ in dominated coloring of $P^{\{u,u_1\}}$ and the superedges $P^{\{u,u_i\}}$, for any $1 \leq i \leq s$, have been colored as a dominated
coloring of $P_{k+1}$ such that for any $i\neq i^\prime , i, i^\prime \in \{1,..., s\}$, (note that $c(y)$ is the color of
vertex $y$ in our coloring)
\begin{equation*}
\left(\bigcup_{j=1}^{k-1} c(x_j^{\{u,u_i\}}\cup c(u_i)\right)\bigcap \left(\bigcup_{j=1}^{k-1} 
c(x_j^{\{u,u_{i^\prime}\}}\cup c(u_{i^\prime})\right)=\varnothing,~ \text{where} ~  i\neq i^\prime.
\end{equation*}
Thus we need at most $(s-1)\chi_{dom}(P_k) + \chi_{dom}(P_{k+1})$ colors for such coloring of vertices of
superedges $P^{\{u,u_i\}}$, $1 \leq i \leq s$. Note that we need at most $\chi_{dom}(P_k)$ new colors for a
dom-coloring of $P^{\{u,u_i\}}$, since the vertex $u$ has been colored in all superedges $P^{\{u,u_i\}}$,
$1 \leq i \leq s$. We do not use the colors used for superedges $P^{\{u,u_i\}}$, $1 \leq i \leq s$, any
more. In the next step, we consider that superedges in $G^{\frac{1}{k}}$ which are replaced instead
of incident edges to $u_i^,s$ in $G$, and have not been colored in the prior step. Now we
color the vertices of these superedges as a dom-coloring of $P_k$, such that the vertices
$u_2,...,u_s$ have been colored in prior step, and the pairwise intersection of the set of
colors used for coloring of vertices of these superedges is the empty set. We continue
this process to color all vertices of $G^{\frac{1}{k}}$. This coloring is a dom-coloring of  $G^{\frac{1}{k}}$, because
every superedge have been colored with distinct color set, except the end vertices of the
superedges, possibly. Finally, since we used at most $(m -1)\chi_{dom}(P_k) + \chi_{dom}(P_{k+1})$ colors,
the right inequality follows.\\
For the left inequality, if $m=1$, then the result is hold. Let $m\geq 2$ and $deg(u) \geq 2$. Let $w, z$ be adjacent with $u$. We consider the superedge $P^{\{u,w\}}$ of $G^{\frac{1}{k}}$. we color this superedge with $\chi_{dom}(P_{k+1})$ colors. Thus we can not use no colors of $P_{k+1}$ for vertex $z$. Therefore we need at least $\chi_{dom}(P_{k+1})$ color for dominated coloring of $G^{\frac{1}{k}}$.\qed
\begin{remark}
 The upper bound in Theorem \ref{frac} is sharp. Consider the graph $K_{1,n}^{\frac{1}{3}}$.
 \end{remark}
 \begin{theorem}\label{dfrac}
 If $G$ is a connected graph with $m$ edges and maximum degree $\Delta (G)$, then
 \begin{equation*}
2+\Delta (G)\chi_{dom}(P_{k-1})-1\leq \chi_{dom}(G^{\frac{1}{k}}) \leq 2+\Delta (G)\chi_{dom}(P_{k-1})-1+(m-\Delta (G))\chi_{dom}(P_{k}).
\end{equation*}
 \end{theorem}
 \proof 
 For the right inequality, let $u$ be the vertex with maximum degree and $N_G(u)=\{u_1, ... , u_{\Delta (G)}\}$. We consider the graph $G^{\frac{1}{k}}$ and color the vertex $u$ with  $1$  and color the adjacent vertices of  $u$ in superedges $P^{\{u,u_i\}}$ ($1 \leq i \leq \Delta (G)$)  with $2$. Now we consider the superedge $P^{\{u,u_1\}}$ with vertices $\{u, x_1^{\{u,u_1\}}, ... , x_{k-1}^{\{u,u_1\}}, u_1 \}$ and put $c(x_2^{\{u,u_1\}})=1$ and the other vertices of superedge $P^{\{u,u_1\}}$ color the same as dominated coloring $P_{k-2}$ with new color. Now remain vertices of other superedges $P^{\{u,u_i\}}$ for $2\leq i \leq \Delta (G)$ color with $(\Delta (G)-1)\chi_{dom}(P_{k-1})$. We do not use the colors used for superedges $P^{\{u,u_i\}}$, $1 \leq i \leq \Delta(G)$, any
more. In the next step, we consider that superedges in $G^{\frac{1}{k}}$ which are replaced instead
of incident edges to $u_i^,s$ in $G$, and have not been colored in the prior step. Now we
color the vertices of these superedges as a dom-coloring of $P_k$, such that the vertices
$u_2, . . . , u_{\Delta(G)}$ have been colored in prior step, and the pairwise intersection of the set of
colors used for coloring of vertices of these superedges is the empty set. We continue
this process to color all vertices of $G^{\frac{1}{k}}$. This coloring is a dom-coloring of  $G^{\frac{1}{k}}$, because
every superedge have been colored with distinct color set, except the end vertices of the
superedges, possibly. Finally, since we used at most $2+\Delta (G)\chi_{dom}(P_{k-1})-1+(m-\Delta (G))\chi_{dom}(P_{k})$ colors,
the right inequality follows. 
The proof of the left inequality is the same of the proof of right inequality, but here we need at least $2+\Delta (G) \chi_{dom}(P_{k-1})-1$ color. \qed
 \begin{remark}
 By considering  the graph $K_{1,n}^{\frac{1}{5}}$, we see that the  bounds in Theorem \ref{dfrac} are sharp. 
 \end{remark}
 
\end{document}